\documentclass[dvips,aos,preprint]{imsart}

\RequirePackage[OT1]{fontenc}
\RequirePackage{amsthm,amsmath}
\RequirePackage[numbers]{natbib}
\RequirePackage[colorlinks,citecolor=blue,urlcolor=blue]{hyperref}
\RequirePackage{hypernat}
\newcommand{\rN}{\mathop{{\rm I}\mskip-4.0mu{\rm
N}}\nolimits}
\arxiv{math.PR/0000000}

\startlocaldefs
\numberwithin{equation}{section}
\theoremstyle{plain}
\newtheorem{thm}{Theorem}[section]
\endlocaldefs

\begin{document}

\begin{frontmatter}
\title{Statistical analysis of the Hirsch index}
\runtitle{Statistical analysis of the Hirsch index}

\begin{aug}
\author{\fnms{Luca} \snm{Pratelli}\thanksref{m1}\ead[label=e1]{luca\_pratelli@marina.difesa.it}},
\author{\fnms{Alberto} \snm{Baccini}\thanksref{m2}\ead[label=e2]{baccini@unisi.it}},
\author{\fnms{Lucio} \snm{Barabesi}\thanksref{m2}\ead[label=e3]{barabesi@unisi.it}}

\and
\author{\fnms{Marzia} \snm{Marcheselli}\thanksref{m2}\ead[label=e4]{marcheselli@unisi.it}}

\runauthor{L. Pratelli et al.}

\affiliation{Accademia Navale di Livorno\thanksmark{m1} and Universit\`a di Siena\thanksmark{m2}}
\address{Address of the First author\\
Accademia Navale, Viale Italia 72,  \\
57100 Livorno, Italy \\
\printead{e1}
\phantom{E-mail:\ }}

\address{Address of the Others authors\\
Dipartimento di Economia Politica,  \\
P.zza S.Francesco 17, 53100 Siena, Italy\\
\printead{e2}\\
\printead{e3}\\
\printead{e4}} 
\end{aug}

\begin{abstract}
The Hirsch index (commonly referred to as $h$-index) is a
bibliometric indicator which is widely recognized as effective for measuring
the scientific production of a scholar since it summarizes size and impact of
the research output. In a formal setting, the $h$-index is actually an
empirical functional of the distribution of the citation counts received by
the scholar.

Under this approach, the asymptotic theory for the empirical
$h$-index has been recently exploited when the citation counts follow a
continuous distribution and, in particular, variance estimation has been
considered for the Pareto-type and the Weibull-type distribution families.
However, in bibliometric applications, citation counts display a distribution
supported by the integers. Thus, we provide general properties for the
empirical $h$-index under the small- and large-sample settings.

In addition,
we also introduce consistent nonparametric variance estimation, which allows
for the implemention of large-sample set estimation for the theoretical
$h$-index.
\end{abstract}

\begin{keyword}[class=AMS]
\kwd[Primary ]{62G05}
\kwd[; secondary ]{62G20, 62G32}
\end{keyword}

\begin{keyword}
\kwd{Hirsch index}
\kwd{heavy-tailed distribution}
\kwd{variance estimation}
\end{keyword}

\end{frontmatter}

\section{Introduction} 
 The $h$-index has been introduced by Hirsch (2005) as a
research performance indicator for individual scholars. The $h$-index is
designed as a single score, balancing the two most important dimensions of
research activity, {\sl i.e.} the productivity of a scholar and the corresponding
impact on the scientific community. Indeed, according to the original
definition of the empirical $h$-index provided by Hirsch (2005), ``a
scientist has index $h$, if $h$ of his or her $N_p$ papers have at least $h$
citations each, whereas the other $(N_p-h)$ papers have no more than $h$
citations each''. 

Notwithstanding that the $h$-index has been only recently proposed, it is
increasingly being adopted for evaluation and comparison purposes to provide
information for funding and tenure decisions, since it is considered an
appropriate tool for identifying ``good'' scientists (Ball, 2007). As a
matter of fact, several reasons explain its popularity and diffusion (Costas
and Bordons, 2007). As it is apparent from its definition, the $h$-index
displays a simple structure allowing easy computation, using data from Web of
Science, Scopus or Google Scholar, while it is robust to publications with a
large or small number of citations. In addition, the $h$-index may be adopted
for assessing the research performance of more complex structures, such as
journals (setting up as a competitor of the Impact Factor, see {\sl e.g.} Braun {\sl et
al.}, 2006), groups of scholars, departments and institutions (Molinari and
Molinari, 2008) and even countries (Nejati and Hosseini Jenab, 2010).

Quite obviously, the $h$-index has received considerable attention by
researchers in the fields of scientometrics and information science (Van
Noorden, 2010). Even if the Hirsch index was originally introduced in a
descriptive framework, scientometricians often aim to assume a statistical
model for citation distribution and the interest is focused on the empirical
$h$-index (see {\sl e.g.} Gl\"anzel, 2006). In a proper statistical perspective,
Beirlant and Einmahl (2010) have managed the empirical $h$-index as the
estimator for a suitable statistical functional of the citation-count
distribution. Accordingly, these authors have proven the consistency of the
empirical $h$-index and they have given the conditions for its large-sample
normality. In addition, Beirlant and Einmahl (2010) have provided an
expression for the large-sample variance of the empirical $h$-index and a
simplified formula for the same quantity when the underlying citation-count
distribution displays Pareto-type or Weibull-type tails. These two special
families have central importance in bibliometrics, since heavy-tailed
citation-count distributions are commonly assumed (see {\sl e.g.} Gl\"anzel, 2006,
and Barcza and Telcs, 2009).

Beirlant and Einmahl (2010) have developed the asymptotic theory for the
empirical $h$-index by assuming a continuous citation-count distribution,
even if the citation number is obviously an integer. Hence, scientometricians
may demand results on the empirical Hirsch index under a more general
approach. Thus, on the basis of a suitable reformulation of the empirical
$h$-index, we provide distributional properties, as well as exact expressions
for the mean and variance, of the empirical $h$-index when citation counts
follow a distribution supported by the integers. Moreover, the general
large-sample properties of the empirical $h$-index are obtained and the link
between the ``integer'' and the ``continuous'' cases is fully analyzed. In
addition, a simple and consistent nonparametric estimator for the variance of
the empirical $h$-index is also introduced under very mild conditions.
Accordingly, the achieved theoretical results are assessed in a small-sample
study by assuming some specific heavy-tailed distributions for the citation
counts. Finally, an application to the ``top-ten'' researchers for the Web of
Science archive in the field of Statistics and Probability during the period
1985-2010 is carried out.

\section{Definitions and preliminary results} Let $X$ be a positive random variable
(r.v.) and let $S$ be the corresponding survival function (s.f.), {\sl i.e.}

\begin{equation*}
S(x)=P(X>x).
\end{equation*}
Even if $X$ is a discrete r.v. in the common bibliometric applications (since
it represents the citation number for a paper of a given scholar), we
actually provide a more general approach. Similarly to Beirlant and Einmahl
(2010), it is assumed that $S(x)>0$ for each $x$ since an
unbounded support for the r.v. $X$  is usually required in scientometrics (Egghe, 2005). If the
left-hand limit of $S$ is denoted by

\begin{equation*}
S_{-}(x)=P(X\geq x),
\end{equation*}
on the basis of the Hirsch definition of the empirical $h$-index reported in
the Introduction, for each integer $n\geq 1$, a ``natural'' expression for
the theoretical Hirsch index $h_n$, corresponding to the law of $X$, is given
by

\begin{equation}
h_n=\sup\{x\geq 0:nS_{-}(x)\geq x\}.
\end{equation}
Obviously, it turns out that $h_n>0$ since $S$ is a strictly positive
function. It is at once apparent that (2.1) encompasses the definition of the
theoretical $h$-index given by Beirlant and Einmahl (2010) when the r.v. $X$
is assumed to be continuous. Moreover, when the r.v. $X$ is integer-valued -
the most interesting situation in bibliometrics - the theoretical Hirsch
index (2.1) reduces to the integer number defined by

\begin{equation}
\begin{aligned}[t]
h_n&=\max\{j\in {\rN}:nS(j-1)\geq j\}\\
&=\sum_{j=1}^nI_{[j/n,1]}(S(j-1)),
\end{aligned}
\end{equation}
where $I_E$ represents the usual indicator function of a set $E$. It should
be remarked that $h_n\nearrow\infty$ and $h_n/n\rightarrow 0$ as
$n\rightarrow\infty$, as immediately follows from the definition (2.1). 

Since it holds that

\begin{equation*}
S_{-}(j)=P(\lfloor X\rfloor\geq j), j\in {\rN},
\end{equation*}
where $\lfloor\cdot\rfloor$ denotes the function giving the greatest integer
less than or equal to the function argument, it is worth noticing that
$\lfloor h_n\rfloor$ turns out to be the $h$-index corresponding to the law
of $\lfloor X\rfloor$. Indeed, from the definition (2.1) we have

\begin{equation*}
S_{-}(\lfloor h_n\rfloor)\geq\frac{\lfloor h_n\rfloor}{n}
\end{equation*}
and

\begin{equation*}
S_{-}(1+\lfloor h_n\rfloor)<\frac{1+\lfloor h_n\rfloor}{n}.
\end{equation*}
Rephrasing the previous statement in its dual setting, if $X$ is an
integer-valued r.v., the $h$-index corresponding to the law of $X$ turns out
to be the integer part of the $h$-index corresponding to the absolutely
continuous law $X+U$, where $U$ is a uniform r.v. on $[0,1]$ independent from
$X$.

If $X_1,\ldots,X_n$ are $n$ independent copies of $X$, the estimator of
$h_n$, {\sl i.e.} the empirical $h$-index, may be immediately introduced as an
empirical functional on the basis of the definition (2.1). More precisely, the
empirical $h$-index is defined to be

\begin{equation}
\widehat{H}_n=\sup\{x\geq 0:n\widehat{S}_{n-}(x)\geq x\},
\end{equation}
where

\begin{equation*}
\widehat{S}_{n-}(x)=\frac{1}{n}\,\sum_{i=1}^nI_{[x,\infty[}(X_i).
\end{equation*}
It should be remarked that (2.3) reduces to the empirical $h$-index defined by
Beirlant and Einmahl (2010) when the r.v. $X$ is continuous. Moreover, on the
contrary to (2.3), the expression of the empirical $h$-index commonly given in
bibliometric literature (see {\sl e.g.} Gl\"anzel, 2006, p.316) is not consistent
when the realizations of the $n$ copies are null. In addition, by considering
the previous discussion and from expression (2.2), the estimator of $\lfloor
h_n\rfloor$ corresponding to the law of $\lfloor X\rfloor$ is given by

\begin{equation}
\widetilde{H}_n=\sum_{j=1}^nI_{[j/n,1]}(\widehat{S}_n(j-1)),
\end{equation}
where $\widehat{S}_n$ represents the empirical s.f., {\sl i.e.}

\begin{equation*}
\widehat{S}_n(x)=\frac{1}{n}\,\sum_{i=1}^nI_{]x,\infty[}(X_i).
\end{equation*}
It is at once apparent that 
\begin{equation}
\widetilde{H}_n=\lfloor\widehat{H}_n\rfloor,
\end{equation}
while it turns out that $\widetilde{H}_n=\widehat{H}_n$ when the r.v. $X$ is
integer-valued. Actually, estimator (2.4) constitutes the formal expression of
the empirical $h$-index given by Hirsch and reported in the Introduction. 

It holds that $\widehat{H}_n\nearrow\infty$ a.s. (and hence
$\widetilde{H}_n\nearrow\infty$ a.s.) as $n\rightarrow\infty$ on the basis of
the Glivenko-Cantelli Theorem. In particular, it follows that
$E[\widehat{H}_n]\nearrow\infty$ and
$E[\widetilde{H}_n]\nearrow\infty$ as $n\rightarrow\infty$.

In order to achieve some useful small-sample properties for the estimator
(2.4), it should be remarked that 
\begin{equation*}
Y_{j,n}=I_{[j/n,1]}(\widehat{S}_n(j-1))\text{ , }j=1,\ldots,n,
\end{equation*}
are dependent Bernoulli random variables. More precisely, each $Y_{j,n}$
turns out to be a Bernoulli r.v. with parameter

\begin{equation*}
\begin{aligned}[t]
p_{j,n}=E[Y_{j,n}]&=P(n\widehat{S}_n(j-1)\geq j)\\
&=\sum_{y=j}^n\binom{n}{y}\,S(j-1)^y(1-S(j-1))^{n-y}.
\end{aligned}
\end{equation*}
In the sequel, we pose $p_{j,n}=0$ if $j>n$. Moreover, since it trivially
holds that
\begin{equation*}
{\rm Var}[Y_{j,n}]=p_{j,n}(1-p_{j,n}),
\end{equation*}
and
\begin{equation*}
{\rm Cov}[Y_{j,n},Y_{l,n}]=p_{j,n}(1-p_{l,n})
\end{equation*}
for $j>l$, it also follows that
\begin{equation*}
E[\widetilde{H}_n]=\sum_{j=1}^np_{j,n}
\end{equation*}
and
\begin{equation}
\begin{aligned}[t]
{\rm Var}[\widetilde{H}_n]&=\sum_{j=1}^np_{j,n}(1-p_{j,n})+2\,\sum_{l=2}^np_{l,n}\,%
\sum_{j=1}^{l-1}(1-p_{j,n})\\
&=\sum_{j=1}^nr_{j,n}(1-p_{j,n}),
\end{aligned}
\end{equation}
where
\begin{equation*}
r_{j,n}=p_{j,n}+2\,\sum_{l=j+1}^np_{l,n}.
\end{equation*}
Obviously, it holds that
$E[\widehat{H}_n]/E[\widetilde{H}_n]\rightarrow 1$ as
$n\rightarrow\infty$. The behavior of ${\rm Var}[\widehat{H}_n]$ and
${\rm Var}[\widetilde{H}_n]$ as $n\rightarrow\infty$ will be considered at length
in Sections 3 and 4.

\section{Large-sample properties of the empirical $h$-index}  By means
of expression (2.5) and considering the discussion following expression (2.2), in
order to explore the large-sample behavior of the empirical $h$-index as
$n\rightarrow\infty$, laws defined on a continuous support may be managed by
considering laws concentrated on integers and {\sl vice versa}. Moreover, if
$(a_n)_n$ is an infinitesimal sequence and $a_n(\widehat{H}_n-h_n)$ converges
in distribution to $\mu$, it follows that

\begin{equation*}
a_n(\widehat{H}_n-h_n)\overset{d}{\longrightarrow}\mu\Longleftrightarrow
a_n(\widetilde{H}_n-\lfloor h_n\rfloor)\overset{d}{\longrightarrow}\mu.
\end{equation*}
In addition, by noting that $a_n\sim b_n$ means asymptotic equivalence of the
sequences $(a_n)_n$ and $(b_n)_n$, {\sl i.e.} $\lim_na_n/b_n=1$ as
$n\rightarrow\infty$, if $a_n^{-2}\sim{\rm Var}[\widetilde{H}_n]$ and
$\lim_n{\rm Var}[\widetilde{H}_n]=\infty$, we have

\begin{equation*}
\lim_n\frac{{\rm Var}[\widehat{H}_n]}{{\rm Var}[\widetilde{H}_n]}=1
\end{equation*}
and

\begin{equation*}
a_n(\widehat{H}_n-h_n)\overset{d}{\longrightarrow}\mu\Longleftrightarrow%
\frac{\widetilde{H}_n-h_n}{\widetilde\sigma_n}\overset{d}{\longrightarrow}\mu%
\Longleftrightarrow\frac{\widehat{H}_n-h_n}{\widetilde\sigma_n}\overset{d}{\longrightarrow}%
\mu,
\end{equation*}
where $\widetilde\sigma_n^2$ is a consistent estimator of ${\rm Var}[\widetilde{H}_n]$, {\sl i.e.}

\begin{equation*}
\frac{\widetilde\sigma_n^2}{{\rm Var}[\widetilde{H}_n]}\overset{P}{\longrightarrow}1.
\end{equation*}

\noindent Hence, in order to implement confidence sets for $h_n$, the evaluation and
the estimation of ${\rm Var}[\widetilde{H}_n]$ is of central importance in the
most interesting case for the scientometricians, {\sl i.e.} when it holds that
${\rm Var}[\widetilde{H}_n]\rightarrow\infty$ as $n\rightarrow\infty$. For
example, this setting happens for the Pareto-type family of laws satisfying
the condition
\begin{equation*}
S(x)=x^{-\alpha}l(x)
\end{equation*}
with $\alpha\in]0,\infty[$ and for the Weibull-type family of laws satisfying
the condition
\begin{equation*}
S(x)=\exp(-x^\tau l(x))
\end{equation*}
with $\tau\in]0,1/2[$, where $l$ is a slowly-varying function, {\sl i.e.}
\begin{equation*}
\frac{l(tx)}{l(t)}\rightarrow 1
\end{equation*}
for each $x$ as $t\rightarrow\infty$. Since the variance (2.6) is a function of the probabilities $p_{j,n}$s, the
preliminary step consists in determining tight inequalities for these
quantities, as the following result provides.

\begin{thm}
If $G$ represents the s.f. of the standard Normal
distribution, there exists a constant $A>0$ such that, for each $n\geq 1$ and
$j=1,\ldots,n$, it holds that

\begin{equation}
|p_{j,n}-G(x_{j,n})|\leq
A\,\frac{v_{j,n}^3+1}{v_{j,n}^4(1+|x_{j,n}|)^6},
\end{equation}
where

\begin{equation*}
x_{j,n}=\frac{j-nS(j-1)}{v_{j,n}}
\end{equation*}
and

\begin{equation*}
v_{j,n}^2=nS(j-1)(1-S(j-1)).
\end{equation*}
\end{thm}
\vskip 0.4cm

 {\scshape Corollary} 3.1.
{\itshape There exists a constant $C>0$ solely depending on $A$, such that
\begin{equation}
\sum_{j=\lfloor 2h_n\rfloor+1}^nr_{j,n}\leq\frac{C}{h_n^{3/2}}
\end{equation}
for each $n$ and $h_n\geq 1$.}
\vskip 0.4cm

The further Corollary to Theorem 3.1 gives the consistency of the estimator
$\widehat{H}_n$.
\vskip0.3cm

{\scshape Corollary} 3.2. {\itshape The ratio $\widehat{H}_n/h_n$ converges in quadratic mean to
$1$, {\sl i.e.} it holds that
\begin{equation*}
E\left[\left(\frac{\widehat{H}_n}{h_n}-1\right)^2\right]%
\rightarrow 0
\end{equation*}
as $n\rightarrow\infty$.}
\vskip 0.3cm

Similarly to the framework considered by Beirlant and Einmahl (2010), the
previous consistency result is stated in a ratio-setting since
$h_n\nearrow\infty$ as $n\rightarrow\infty$. Finally, on the basis of
Corollary 3.2  it also follows that
\begin{equation*}
\frac{E[\widehat{H}_n]}{h_n}\rightarrow 1
\end{equation*}
as $n\rightarrow\infty$. 

\section{Consistent estimation of the empirical
$h$-index variance}  As emphasized in Section 2, in order to
achieve the convergence in distribution of $\widehat{H}_n$, the evaluation of
the variance (2.6) is central. To this aim, the following result provides some
inequalities and a useful asymptotic equivalence for (2.6) by assuming a mild
condition.

\begin{thm}
 For each $n$ it holds that
\begin{equation}
{\rm Var}[\widetilde{H}_n]\geq\sum_{j=1}^{\lfloor
2h_n\rfloor}r_{j,n}(1-p_{j,n})\geq\sum_{j=1}^{\lfloor
2h_n\rfloor}\widetilde{r}_{j,n}(1-p_{j,n}),
\end{equation}
where
\begin{equation*}
\widetilde{r}_{j,n}=p_{j,n}+2\,\sum_{l=j+1}^{\lfloor
2h_n\rfloor}p_{l,n}.
\end{equation*}
Moreover, if
\begin{equation}
\liminf_n{\rm Var}[\widetilde{H}_n]>0,
\end{equation}
it holds that
\begin{equation}
{\rm Var}[\widetilde{H}_n]\sim\sum_{j=1}^{\lfloor
2h_n\rfloor}\widetilde{r}_{j,n}(1-p_{j,n})
\end{equation}
as $n\rightarrow\infty$. In particular, if
\begin{equation*}
V_n=\sum_{j=1}^{\lfloor
3\widehat{H}_n\rfloor}p_{j,n}(1-p_{j,n})+2\,\sum_{l=2}^{\lfloor
3\widehat{H}_n\rfloor}p_{l,n}\,\sum_{j=1}^{l-1}(1-p_{j,n})
\end{equation*}
it holds that
\begin{equation*}
R_n=\frac{V_n}{{\rm Var}[\widetilde{H}_n]}\overset{P}{\longrightarrow}1
\end{equation*}
as $n\rightarrow\infty$.
\end{thm}

From Teorem 4.1, it is at once apparent that $V_n$ would be a consistent
estimator of (2.6) when it is possible to evaluate the $p_{j,n}$s for
$j\leq\lfloor 3\widehat{H}_n\rfloor$ and in the case that condition (4.2)
holds, {\sl i.e.} when ${\rm Var}[\widetilde{H}_n]$ does not approach $0$ as
$n\rightarrow\infty$. Hence, it is convenient to introduce a further
condition which solely involves the behavior of $S$ on $\rN$ and which
implies condition (4.2). More precisely, we consider the condition

\begin{equation}
\lim_n\frac{\sqrt{n}\psi(n)}{S(n)}=0,
\end{equation}
where
\begin{equation*}
\psi(n)=S(n-1)-S(n)=P(n-1<X\leq n).
\end{equation*}
Obviously, when the r.v. $X$ is integer-valued, $\psi$ represents the
probability function corresponding to $X$. It should be noticed that
condition (4.4) may also be reformulated as

\begin{equation*}
\frac{S(n-1)}{S(n)}=1+\frac{\gamma_n}{\sqrt{n}},
\end{equation*}
where $(\gamma_n)_n$ is a positive infinitesimal sequence, and hence for each
$M>0$, it holds
\begin{equation}
\lim_n\frac{S(n-M\sqrt{n})}{S(n+M\sqrt{n})}=1,
\end{equation}
since

\begin{equation*}
\begin{aligned}[t]
1&\leq\frac{S(n-M\sqrt{n})}{S(n+M\sqrt{n})}\sim\prod_{j=\lfloor-M\sqrt{n}%
\rfloor+1}^{\lfloor
M\sqrt{n}\rfloor+1}\left(1+\frac{\gamma_{n+j}}{\sqrt{n+j}}\right)\\
&\leq\exp\left(\frac{2(M+1)\sqrt{n}\delta_{n+\lfloor-M\sqrt{n}\rfloor}}{%
\sqrt{n-M\sqrt{n}-1}}\right)\sim 1
\end{aligned}
\end{equation*}
where $\delta_n=\sup_{h\geq n}\gamma_h$. As proven in the following result,
condition (4.4) ensures the unboundedness of (2.6) as $n\rightarrow\infty$.

\begin{thm}
If the law of $X$ satisfies condition (4.4), it holds that

\begin{equation*}
\lim_n{\rm Var}[\widetilde{H}_n]=\infty,
\end{equation*}
\end{thm}
\ \
\noindent In order to achieve consistent estimation of (2.6), it is necessary to
introduce a further condition, which is slightly more restrictive than
condition (4.4). More precisely, this condition assumes that for each $M\geq
0$ it holds that
\begin{equation}
\lim_n\,\left(\sup_{j\in
D_{M,n}}\left|\frac{\psi(j)}{\psi(n)}-1\right|\right)=0,
\end{equation}
where $D_{M,n}=[n-M\sqrt{n},n+M\sqrt{n}]\cap\rN$. It may be easily
verified that condition (4.6) implies condition (4.4) and hence condition (4.5).

Since a natural estimator for $p_{j,n}$ is given by
\begin{equation*}
\widehat{p}_{j,n}=\sum_{y=j}^n\binom{n}{y}\,\widehat{S}_n(j-1)^y(1-%
\widehat{S}_n(j-1))^{n-y},
\end{equation*}
on the basis of the large-sample behavior of the ratio $R_n$ given in
Theorem 4.1, an estimator for the variance (2.6) turns out to be
\begin{equation}
\widehat{V}_n=\sum_{j=1}^{\min(\lfloor
3\widehat{H}_n\rfloor,n)}\widehat{p}_{j,n}(1-\widehat{p}_{j,n})+2\,%
\sum_{l=2}^{\min(\lfloor
3\widehat{H}_n\rfloor,n)}\widehat{p}_{l,n}\,\sum_{j=1}^{l-1}(1-\widehat{p}_{j,%
n}).
\end{equation}
It should be remarked that estimator (4.7) is fully nonparametric. Indeed, it
does not require the specification of a semi-parametric model for the
underlying citation distribution as in the case of the variance estimator
proposed by Beirlant and Einmahl (2010). For example, their estimator
requires the estimation of the Paretian index when a Pareto-type family is
assumed for the law of $X$ - a non-trivial task, see {\sl e.g.} Beirlant {\sl et al.}
(2004). 

The following result provides a compact asymptotic equivalent
expression for (2.6) and the consistency of estimator (4.7) if condition (4.6) is
assumed.

\begin{thm}
 If the law of $X$ satisfies condition (4.6), it holds that
\begin{equation}
{\rm Var}[\widehat{H}_n]\sim\frac{h_n}{(1+n\psi(\lfloor h_n\rfloor))^2}
\end{equation}
as $n\rightarrow\infty$. Moreover, it follows that
\begin{equation*}
\frac{\widehat{V}_n}{{\rm Var}[\widehat{H}_n]}\overset{P}{\longrightarrow}1
\end{equation*}
as $n\rightarrow\infty$.
\end{thm}

It should be remarked that for the Pareto-type and the Weibull-type families
(described in Section 3) condition (4.6) is satisfied. Accordingly,
$\widehat{H}_n$ approaches normality and from Theorem 4.3 for the
Pareto-type family, it holds that
\begin{equation*}
{\rm Var}[\widehat{H}_n]\sim\frac{h_n}{(1+\alpha)^2}
\end{equation*}
and
\begin{equation*}
h_n\sim Cn^{1/(1+\alpha)}
\end{equation*}
when $l(x)\sim C^{1+\alpha}$, while for the Weibull-type family, it holds that
\begin{equation*}
{\rm Var}[\widehat{H}_n]\sim\frac{h_n}{(1+\tau\log(n/h_n))^2}
\end{equation*}
and
\begin{equation*}
h_n\sim C(\log n)^{1/\tau}\text{ ,}
\end{equation*}
when $l(x)\sim C^{-\tau}$ and where $C>0$ is a suitable constant. These
results are in complete agreement with the findings by Beirlant and Einmahl
(2010).

On the basis of Theorem 4.3 and on the remarks contained in Section 2, when
condition (4.6) is satisfied by the underlying distribution, a large-sample
confidence set for $h_n$ at the $(1-\gamma)$ confidence level turns out to be
\begin{equation}
C_n=\{[\![\widehat{H}_n-z_{1-\gamma/2}\sqrt{\widehat{V}_n}]\!],\ldots,[\![%
\widehat{H}_n+z_{1-\gamma/2}\sqrt{\widehat{V}_n}]\!]\},
\end{equation}
where $z_\gamma$ represents the $\gamma$-th quantile of the standard Normal
distribution, while $[\![\cdot]\!]$ represents the function giving the
integer closest to the argument. In addition, in order to assess the
homogeneity of the theoretical $h$-indexes for two scholars, a suitable test
statistic is given by
\begin{equation*}
T_n=\frac{\widehat{H}_{1,n}-\widehat{H}_{2,n}}{\sqrt{\widehat{V}_{1,n}+%
\widehat{V}_{2,n}}},
\end{equation*}
where $\widehat{H}_{1,n}$ and $\widehat{H}_{2,n}$ represent the empirical
$h$-indexes corresponding to the scholars, while $\widehat{V}_{1,n}$ and
$\widehat{V}_{2,n}$ are the respective variance estimators as given by (4.7).
It is at once apparent that 
\begin{equation*}
T_n\overset{d}{\longrightarrow}{\mathcal N}(0,1)
\end{equation*}
as $n\rightarrow\infty$, when $\widehat{H}_{1,n}$ and $\widehat{H}_{2,n}$
approaches normality. The test statistic $T_n$ is defined in a nonparametric
setting, in contrast to the test statistic proposed in a semiparametric
approach by Beirlant and Einmahl (2010), which requires consistent estimation
of the two Paretian indexes of the scholar citation distributions. 

Finally, when the analysis of the theoretical $h$-indexes corresponding to $k$
scholars ($k\geq 2$) is considered, simultaneous set estimation and
homogeneity hypothesis testing could be managed by means of techniques
similar to those suggested in Marcheselli (2003). These issues will be
pursued in future research.

\section{Some studies and examples} In order to assess in practice the properties
of the empirical $h$-index achieved in the previous sections, a study was
carried out for two heavy-tailed distributions. First, a discrete stable
distribution for the r.v. $X$ was considered. This distribution may be
specified via the probability generating function
\begin{equation*}
g(s)=\text{E}[s^X]=\exp(-\lambda(1-s)^\alpha)\text{ , }s\in[0,1],
\end{equation*}
where $\alpha\kern-1mm\in\, ]0,1]$ and $\lambda\kern-1mm\in\,]0,\infty[$ (Steutel and van Harn,
2004, p.265). The discre\-te stable distribution is Paretian of order $\alpha$
for $\alpha\kern-1mm\in\,]0,1[$ (Christoph and Schreiber, 1998) and it constitutes a
flexible and natural model for heavy-tailed discrete data (see Marcheselli {\sl et
al.}, 2008, for a description of the distribution properties and of the
corresponding parameter estimation issues). A ``discretized'' Weibull
distribution was subsequently assumed for the r.v. $X$. The distribution
displays the probability function
\begin{equation*}
f(x)=[\exp(-x^\tau)-\exp(-(x+1)^\tau)]I_{\rN}(x),
\end{equation*}
where $\tau\in]0,\infty[$. Obviously, it turns out that
\begin{equation*}
S(j)=\exp(-(j+1)^\tau), \ j\in\rN.
\end{equation*}
By assuming that $n=30,50,100,150,200$, the values of $h_n$,
$E[\widehat{H}_n]$, ${\rm Var}[\widehat{H}_n]$ and
of the large-sample variance approximation $(4.8)$ were computed for the discrete stable distribution
with parameter vectors $(\alpha,\lambda)=(0.25,1.0)$, $(0.50,1.5),(0.75,2.0)$,
as well as for the discretized Weibull distribution with parameters
$\tau=0.01,0.10,0.40$. These choices were made in order to fit, as close as
possible, the real productivity and the real citation distributions of
scholars with different scientific ages and belonging to different research
areas and with more or less pronounced impact on research. 

In the study,
$B=10,000$ random variates were generated for each $n$ choice and for each
considered distribution in order to achieve the Monte Carlo estimation of
$E[\widehat{V}_n]$ and the Monte Carlo estimation of the actual
coverage for the confidence set (4.9) at the $95$\% nominal confidence level.
The corresponding results were reported in Tables I and II. 

The analysis of
these tables shows that $h_n$ and $E[\widehat{H}_n]$ are similar
even for small $n$ values and $\widehat{V}_n$ turns out to be nearly
unbiased. In addition, it can be verified that the actual coverage of the
large-sample confidence set (4.9) is almost equivalent to the nominal coverage
even for small $n$ values. Unfortunately, it can be assessed that the
large-sample variance approximation (4.8) may be quite dissimilar from
${\rm Var}[\widehat{H}_n]$ even for quite large $n$ values. It should be remarked
that an estimation procedure based on (4.8) requires, in any case, the
additional estimation of $\alpha$ or $\tau$.

Accordingly, $\widehat{V}_n$
seems to be an appealing variance estimator, both from a theoretical and
practical perspective. In general, we have verified similar conclusions for a
plethora of distributions satisfying condition (4.6), even if we have not
reported the corresponding results.

\vskip0,3cm
\begin{tabular}[t]{ccccccccc}
\multicolumn{9}{c}{{\bf Table I.} Discrete stable distribution.}\\
\hline
$\alpha$&$\lambda$&$n$&$h_n$&$E[\widehat{H}_n]$&${\rm Var}[%
\widehat{H}_n]$&$\frac{h_n}{(1+\alpha)^2}$&$E[%
\widehat{V}_n]$&Coverage\\
\hline
$0.25$&$1.0$&$30$&$11$&$11.31$&$4.73$&$7.04$&$4.88$&$0.96$\\
&&$50$&$17$&$17.17$&$7.58$&$10.88$&$7.79$&$0.96$\\
&&$100$&$30$&$30.28$&$14.18$&$19.20$&$14.51$&$0.96$\\
&&$150$&$42$&$42.19$&$20.34$&$26.88$&$20.73$&$0.95$\\
&&$200$&$53$&$53.38$&$26.21$&$33.92$&$26.74$&$0.95$\\
&&&&&&&&\\
$0.50$&$1.5$&$30$&$9$&$8.96$&$2.66$&$4.00$&$2.96$&$0.95$\\
&&$50$&$12$&$12.46$&$4.01$&$5.33$&$4.35$&$0.96$\\
&&$100$&$19$&$19.59$&$6.83$&$8.44$&$7.25$&$0.96$\\
&&$150$&$25$&$25.57$&$9.25$&$11.11$&$9.72$&$0.96$\\
&&$200$&$31$&$30.91$&$11.43$&$13.78$&$11.98$&$0.95$\\
&&&&&&&&\\
$0.75$&$2.0$&$30$&$6$&$6.65$&$1.11$&$1.96$&$1.38$&$0.97$\\
&&$50$&$8$&$8.38$&$1.58$&$2.61$&$1.90$&$0.96$\\
&&$100$&$11$&$11.67$&$2.57$&$3.59$&$2.93$&$0.96$\\
&&$150$&$14$&$14.29$&$3.38$&$4.57$&$3.76$&$0.95$\\
&&$200$&$16$&$16.54$&$4.09$&$5.22$&$4.54$&$0.95$\\
\hline
\end{tabular}%
\vskip 0,6cm
\begin{tabular}[t]{cccccccc}
\multicolumn{8}{c}{{\bf Table II.} Discretized Weibull distribution.}\\
\hline
$\tau$&$n$&$h_n$&$E[\widehat{H}_n]$&${\rm Var}[\widehat{H}_n]$&$%
\frac{h_n}{(1+\tau\log(n/h_n))^2}$&$E[\widehat{V}_n]$&Coverage\\
\hline
$0.01$&$30$&$10$&$10.77$&$6.77$&$9.78$&$6.56$&$0.94$\\
&$50$&$17$&$17.86$&$11.25$&$16.64$&$11.05$&$0.96$\\
&$100$&$35$&$35.47$&$22.43$&$34.28$&$22.25$&$0.96$\\
&$150$&$52$&$52.98$&$33.57$&$50.92$&$33.39$&$0.96$\\
&$200$&$70$&$70.44$&$44.70$&$68.55$&$44.52$&$0.95$\\
&&&&&&&\\
$0.10$&$30$&$8$&$8.63$&$4.93$&$6.24$&$4.90$&$0.95$\\
&$50$&$13$&$13.60$&$7.83$&$10.10$&$7.82$&$0.95$\\
&$100$&$25$&$25.09$&$14.59$&$19.28$&$14.67$&$0.95$\\
&$150$&$35$&$35.83$&$20.95$&$26.67$&$21.12$&$0.95$\\
&$200$&$46$&$46.07$&$27.05$&$34.97$&$27.20$&$0.95$\\
&&&&&&&\\
$0.40$&$30$&$4$&$4.47$&$1.40$&$1.23$&$1.49$&$0.97$\\
&$50$&$6$&$6.01$&$1.72$&$1.76$&$1.85$&$0.96$\\
&$100$&$8$&$8.74$&$2.22$&$1.98$&$2.38$&$0.96$\\
&$150$&$11$&$10.75$&$2.54$&$2.63$&$2.71$&$0.96$\\
&$200$&$12$&$12.38$&$2.77$&$2.66$&$2.96$&$0.96$\\
\hline
\end{tabular}%
\vskip 0,4cm
As a practical application of the achieved results, we also considered the
scientific performance of the best ten scholars in the field of Statistics
and Probability according to the Web of Science archive. Data were drawn from
the Thomson-Reuters databases by selecting the scholars listed in the
category Mathematics of the ISIHighlyCited.com database (given at the WEB
site http://hcr3.isiknowledge.com/home.cgi). For each scholar in the
database, an author search was performed during the month of December 2010 on
the ISI Web of Science for the period 1985-2010. The search was carried out
in such a way that only articles and letters published in journals contained
in the Statistics and Probability database were considered. Accordingly, the
citation counts were collected for each scholar. The citation counts covered
documents contained in the Science Citation Index Expanded and Social
Sciences Citation Index and Arts \& Humanities Citation Index. Finally, the
ten scholars with the highest $h$-indexes were considered. More precisely,
the names of the ten scholars, the corresponding paper number and $h$-index,
as well as the large-sample confidence sets at the $95$\% nominal confidence
level were reported in Table III. Obviously, pratictioners may largely
benifit from this example in order to understand the importance of
quantifying variability for an appropriate comparison analysis of the
research performance.
\vskip 0,4cm

\begin{tabular}[t]{lccc}
\multicolumn{4}{c}{{\bf Table III.} Performance of the ``top-ten'' most-cited
scholars}\\
\multicolumn{4}{c}{ in the field of Statistics and Probability during the
period 1985-2010.}\\
\hline
&$n$&$h_n$&$C_n$\\
\hline
Hall, P.G.&$418$&$46$&$\{42,\ldots,50\}$\\
Rubin, D.B.&$104$&$39$&$\{32,\ldots,46\}$\\
Carroll, R.J.&$198$&$38$&$\{33,\ldots,43\}$\\
Tibshirani, R.&$104$&$37$&$\{31,\ldots,43\}$\\
Fan, J.&$114$&$36$&$\{30,\ldots,42\}$\\
Marron, J.S.&$107$&$36$&$\{31,\ldots,41\}$\\
Hastie, T.J.&$77$&$34$&$\{27,\ldots,41\}$\\
Lin, D.Y.&$93$&$32$&$\{26,\ldots,38\}$\\
Raftery, A.E.&$88$&$31$&$\{25,\ldots,37\}$\\
Wei, L.J.&$88$&$31$&$\{26,\ldots,36\}$\\
\hline
\end{tabular}%

\vfill\eject

\appendix
\vskip9mm
\section{}\label{app}

\subsection{Proof of Theorem 3.1} For fixed $j$ and $n$, let us assume that

\begin{equation*}
Z_i=\frac{I_{]j-1,\infty[}(X_i)-S(j-1)}{\sqrt{S(j-1)(1-S(j-1))}},\ \ i=1,\ldots,n.
\end{equation*}
Accordingly, we have

\begin{equation*}
p_{j,n}=P\left(\frac{1}{\sqrt{n}}\,\sum_{i=1}^nZ_i\geq x_{j,n}\right).
\end{equation*}
Thus, by applying the Osipov inequality (see {\sl e.g.} DasGupta, 2008, p.659) to
the $Z_i$'s for $\alpha=6$, inequality (3.1) is proven, since for each $m\geq
2$ it holds that

\begin{equation*}
E[|Z_i|^m]\leq\frac{1}{[S(j-1)(1-S(j-1))]^{(m-1)/2}}.\eqno \qedsymbol
\end{equation*}

\subsection{Proof of Corollary 3.1} Since for each $n$ and $j\geq 2h_n+1$, it holds that
\begin{equation}
v_{j,n}x_{j,n}=j\left(1-\frac{nS(j-1)}{j}\right)\geq
j\left(1-\frac{nS(2h_n)}{2h_n}\right)\geq\frac{j}{2}
\end{equation}
and since from the definition of $h_n$ it follows that
\begin{equation}
v_{j,n}^2\leq nS(j-1)\leq j
\end{equation}
for $j\geq 2h_n+1$, by means of inequality (3.1) we have
\begin{equation}
|p_{j,n}-G(x_{j,n})|\leq
A\,\frac{v_{j,n}^5+v_{j,n}^2}{v_{j,n}^6x_{j,n}^6}\leq
64A\,\frac{v_{j,n}^5+v_{j,n}^2}{j^6}\leq\frac{128A}{j^{7/2}}.
\end{equation}
Since on the basis of (A.1) and (A.2) it also holds that

\begin{equation*}
x_{j,n}\geq\frac{\sqrt{j}}{2},
\end{equation*}
and hence $G(x_{j,n})\leq G(\sqrt{j}/2)$, for each $l>2h_n$ it follows from
(A.3) 

\begin{equation*}
\sum_{l=j+1}^np_{l,n}\leq\frac{256A}{5j^{5/2}}+\sum_{l=j+1}^nG(\sqrt{l}/2)%
\leq\frac{B}{j^{5/2}},
\end{equation*}
where $B$ is a constant, solely depending on $A$. Thus, it also turns out that

\begin{equation*}
\sum_{j=\lfloor
2h_n\rfloor+1}^nr_{j,n}\leq\frac{B}{h_n^{5/2}}+\frac{4B}{3h_n^{3/2}}
\end{equation*}
and hence inequality (3.2) follows.\hfill{\qedsymbol}

\subsection{Proof of Corollary 3.2} If $S$ is a continuous function, from Corollary 1 given
by Beirlant and Einmahl (2010) it holds that

\begin{equation*}
\frac{\widehat{H}_n}{h_n}\overset{P}{\rightarrow}1
\end{equation*}
as $n\rightarrow\infty$. When $S$ is not a continuous function, let $h_n'$ be
the theoretical $h$-index corresponding to the law of $\lfloor X\rfloor+U$,
where $U$ is a uniform r.v. on $[0,1]$ independent from $X$, while let
$\widehat{H}_n^{\text{ }\prime}$ be the empirical $h$-index based on $n$
independent copies of $\lfloor X\rfloor+U$. 

Since $|h_n-h_n'|\leq 1$ and
$|\widehat{H}_n-\widehat{H}_n^{\text{ }\prime}|\leq 1$, the convergence in
probability to $1$ of $\widehat{H}_n/h_n$ is in turn obtained from the
continuous-setting result. Moreover, the uniform integrability of
$\widehat{H}_n^2/h_n^2$ follows by considering inequality (3.2) since

\begin{equation*}
\begin{aligned}[t]
E\left[\left(\sum_{j=\lfloor
2h_n\rfloor+1}^nY_{j,n}\right)^2\right]&={\rm Var}\left[\sum_{j=\lfloor
2h_n\rfloor+1}^nY_{j,n}\right]+E\left[\sum_{j=\lfloor
2h_n\rfloor+1}^nY_{j,n}\right]^2\\
&\leq\sum_{j=\lfloor 2h_n\rfloor+1}^nr_{j,n}+\left(\sum_{j=\lfloor
2h_n\rfloor+1}^nr_{j,n}\right)^2
\end{aligned}
\end{equation*}
and

\begin{equation*}
\sum_{j=1}^{\lfloor 2h_n\rfloor}Y_{j,n}\leq 2h_n.\eqno{\qedsymbol}
\end{equation*}

\subsection{Proof of Theorem 4.1} The inequalities in (4.1) easily follows from
expression (2.6), while (4.3) follows from (3.2) and

\begin{equation*}
\sum_{j=1}^{\lfloor 2h_n\rfloor}|r_{j,n}-\widetilde{r}_{j,n}|\leq
2h_n\sum_{l=\lfloor 2h_n\rfloor+1}^nr_{l,n}.
\end{equation*}
Moreover, on the basis of Corollary 3.2 it turns out that $\widehat{H}_n/n$
converges in mean to $0$. Hence, it is convenient to consider

\begin{equation*}
R_n'=I_{[0,n]}(3\widehat{H}_n)R_n,
\end{equation*}
for which it holds that $R_n'\leq 1$ by means of (2.6). Hence, the second part
of the Theorem follows from (4.1), (4.3) and

\begin{equation*}
\lim_nP(\lfloor 3\widehat{H}_n\rfloor\leq\lfloor 2h_n\rfloor)=0.\eqno{\qedsymbol}
\end{equation*}

\subsection{Proof of Theorem 4.2} By using the notations introduced in Theorem 4.1,
we have

\begin{equation*}
x_{j+1,n}-x_{j,n}=\frac{1+n\psi(j)}{v_{j+1,n}}+x_{j,n}\,\frac{v_{j,n}-v_{j+1,%
n}}{v_{j+1,n}}.
\end{equation*}
Thus, for a given $c>0$ and for each $n$ such that $h_n\geq 1$, on the basis
of (4.5) it holds that

\begin{equation*}
\frac{1+n\psi(j)}{v_{j+1,n}}<\frac{1}{2c}
\end{equation*}
for $j\in[h_n-\sqrt{h_n},h_n+\sqrt{h_n}]\cap\rN$. Equivalently, there
exist at least $c$ values $x_{j,n}$ in the interval $[-1,1]$. 

Thus, if $D_n$
represents the set of indexes
$j\in[h_n-\sqrt{h_n},h_n+\sqrt{h_n}]\cap\rN$ for which $|x_{j,n}|\leq
1$, on the basis of (3.1) and (4.1) it follows that

\begin{equation*}
\begin{aligned}[t]
{\rm Var}[\widehat{H}_n]&\geq\sum_{j=1}^{\lfloor
2h_n\rfloor}p_{j,n}(1-p_{j,n})\geq\sum_{j\in D_n}p_{j,n}(1-p_{j,n})\\
&\geq cG(1)(1-G(1))-A\,\sum_{j\in
D_n}\frac{v_{j,n}^3+1}{v_{j,n}^4(1+|x_{j,n}|)^6},
\end{aligned}
\end{equation*}
From (4.5) we have

\begin{equation*}
\inf_{j\in
D_n}v_{j,n}\sim\sqrt{nS(h_n-\sqrt{h_n})}\sim\sqrt{nS(h_n)}\sim\sqrt{h_n}%
\text{ ,}
\end{equation*}
and, since $c$ is arbitrary, it holds that

\begin{equation*}
\begin{aligned}[t]
\liminf_n{\rm Var}[\widetilde{H}_n]&\geq
cG(1)(1-G(1))-4A\limsup_n\frac{\sqrt{h_n}}{\inf_{j\in D_n}v_{j,n}}\\
&=cG(1)(1-G(1))-4A,
\end{aligned}
\end{equation*}
which completes the proof.\hfill{\qedsymbol}

\subsection{Proof of Theorem 4.3} For a fixed $M>0$, from condition (4.6) and from (4.5)
it follows that
\begin{equation}
\lim_n\left(\sup_{j\in
D_{M,h_n}}\left|\frac{(1+n\psi(j))v_{h_n,n}}{(1+n\psi(\lfloor
h_n\rfloor))v_{j+1,n}}-1\right|\right)=0.
\end{equation}
Thus, by means of expression (A.4), from Theorem 4.2 it follows that
\begin{equation*}
\begin{aligned}[t]
{\rm Var}[\widehat{H}_n]&\sim{\rm Var}[\widetilde{H}_n]\sim\sum_{j=1}^{\lfloor
2h_n\rfloor}p_{j,n}(1-p_{j,n})+2\,\sum_{l=2}^{\lfloor
2h_n\rfloor}\,p_{l,n}\,\sum_{j=1}^{l-1}(1-p_{j,n})\\
&\sim 2\,\sum_{l=2}^{\lfloor
2h_n\rfloor}\,p_{l,n}\,\sum_{j=1}^{l-1}(1-p_{j,n})\\
&\sim\frac{2h_n}{(1+n\psi(\lfloor
h_n\rfloor))^2}\,\int_{-\infty}^{2h_n}G(x)\,\text{d}x\,\int_{-%
\infty}^x(1-G(u))\,\text{d}u\\
&\sim\frac{h_n}{(1+n\psi(\lfloor
h_n\rfloor))^2},
\end{aligned}
\end{equation*}
since

\begin{equation*}
\int_{-\infty}^\infty
G(x)\,\text{d}x\,\int_{-\infty}^x(1-G(u))\,\text{d}u=\frac{1}{2}.
\end{equation*}
Hence, expression (4.8) is proven. As to the consistency of $\widehat{V}_n$,
by assuming that

\begin{equation*}
\widehat{\psi}_n(j)=\widehat{S}_n(j-1)-\widehat{S}_n(j),
\end{equation*}
it holds that

\begin{equation*}
\sup_{j\in
D_{M,h_n}}\left|\frac{(1+n\widehat{\psi}_n(j))v_{h_n,n}}{(1+n\psi(\lfloor
h_n\rfloor))v_{j+1,n}}-1\right|\overset{P}{\longrightarrow}0,
\end{equation*}
as $n\rightarrow\infty$, since

\begin{equation*}
\inf_{M\geq 0}\limsup_nP(|\widehat{H}_n-h_n|>M\sqrt{h_n})=0.
\end{equation*}
Thus, convergence in probability of $\widehat{V}_n/{\rm Var}[\widetilde{H}_n]$ to
$1$ follows.\hfill{\qedsymbol}

\end{document}